\documentclass[11pt,amsfonts]{article}
\usepackage{graphicx}
\usepackage{latexsym}
\usepackage{amssymb}
\usepackage{amsmath}
\usepackage{enumerate}
\usepackage{layout}
\usepackage{eufrak}
\newtheorem{prop}{Proposition}
\newtheorem{lemma}{Lemma}

\newtheorem{theorem}{Theorem}
\newtheorem{remark}{Remark}

\def\real{{\mathord{{\rm I\kern-2.8pt R}}}}        
\def\inte{{\mathord{{\rm I\kern-2.8pt N}}}}

\def\sZZ{{\rm Z\kern-2.8ptem{}Z}}

\def\z{{\mathchoice
  {\sZZ}
  {\sZZ}
  {\rm Z\kern-0.30em{}Z}
  {\rm Z\kern-0.25em{}Z} }}
\def\sQQ{{\kern 0.27em \vrule height1.45ex width0.03em depth0em
          \kern-0.30em \rm Q}}
\def\qu{{\mathchoice
    {\sQQ}
    {\sQQ}
  {\kern 0.225em \vrule height1.05ex width0.025em depth0em \kern-0.25em \rm Q}
  {\kern 0.180em \vrule height0.78ex width0.020em depth0em \kern-0.20em \rm Q}
        }}
\def\sCC{{\kern 0.27em \vrule height1.45ex width0.03em depth0em
          \kern-0.30em \rm C}}
\def\complex{{\mathchoice
    {\sCC}
    {\sCC}
  {\kern 0.225em \vrule height1.05ex width0.025em depth0em \kern-0.25em \rm C}
  {\kern 0.180em \vrule height0.78ex width0.020em depth0em \kern-0.20em \rm C}
        }}

\newcommand{\ba}{\begin{array}}
\newcommand{\ea}{\end{array}}
\newcommand{\be}{\begin{equation}}
\newcommand{\ee}{\end{equation}}
\newcommand{\bea}{\begin{eqnarray}}
\newcommand{\eea}{\end{eqnarray}}
\newcommand{\beaa}{\begin{eqnarray*}}
\newcommand{\eeaa}{\end{eqnarray*}}

\newcommand{\eps}{\varepsilon}

\def\z{\zeta}

\font\tenmath=msbm10 \font\sevenmath=msbm7 \font\fivemath=msbm5
\newfam\mathfam \textfont\mathfam=\tenmath
\scriptfont\mathfam=\sevenmath \scriptscriptfont\mathfam=\fivemath

\def \={{\buildrel {\rm (law)} \over =}}

\def\qed{ \hfill \vrule width.25cm height.25cm depth0cm\smallskip}

\newcommand{\basa}{\begin{assumption}}
\newcommand{\easa}{\end{assumption}}

\newcommand{\bas}{\begin{assum}}
\newcommand{\eas}{\end{assum}}

\def\limsup{\mathop{\overline{\rm lim}}}

\def\bE{{\bf E}}

\newcommand{\ignore}[1]{}
\allowdisplaybreaks
\begin{document}

\renewcommand{\thefootnote}{\fnsymbol{footnote}}

\title{A strong convergence to the Rosenblatt process \footnote{Dedicated to the memory of Constantin Tudor}}

\author{Johanna Garz\'on $^{1}$ \qquad  Soledad Torres  $^{1,}$\thanks{The author is supported by Dipuv grant 26/2009, ECOS/CONICYT C10E03 2010, MATHAMSUD 09/05/SAMP Stochastic Analysis and Mathematical Physics Research Network.} \qquad
Ciprian A. Tudor $^{2,}$ \footnote{  Partially supported by the ANR grant "Masterie" BLAN 012103. Associate member of the team Samos, Universit\'e de Panth\'eon-Sorbonne Paris 1 }\vspace*{0.1in} \\
$^{1}$ Departamento de Estad\'istica, CIMFAV  Universidad de Valpara\'iso,\\ Casilla
 123-V, 4059 Valparaiso, Chile.
\\margaret.garzon@uv.cl, \vspace*{0.1in} soledad.torres@uv.cl\vspace*{0.1in} \\
 $^{2}$ Laboratoire Paul Painlev\'e, Universit\'e de Lille 1\\
 F-59655 Villeneuve d'Ascq, France.\\
 \quad tudor@math.univ-lille1.fr\vspace*{0.1in}}

 \maketitle

\begin{abstract}
We give a strong approximation of  Rosenblatt process via transport processes and we give the  rate of convergence.
\end{abstract}

\vskip0.5cm

{\bf  2010 AMS Classification Numbers:} 60B10, 60F05, 60H05.

 \vskip0.3cm


\section{Introduction}
Self-similar stochastic processes are of practical interest in various
applications, including econometrics, internet traffic, and hydrology. These
are processes $X=\left(  X\left(  t\right)  :t\geq0 \right)  $ whose
dependence on the time parameter $t$ is self-similar, in the sense that there
exists a (self-similarity) parameter $H\in(0,1)$ such that for any constant
$c>0$, $\left(  X\left(  ct\right)  :t\geq0\right)  $ and $\left(
c^{H}X\left(  t\right)  :t\geq0\right)  $ have the same finite dimensional distributions. These
processes are often endowed with other distinctive properties.

The fractional Brownian motion (fBm) is the usual candidate to model phenomena
in which the self-similarity property can be observed from the empirical data.
This fBm $B^{H}$ is the continuous centered Gaussian process with covariance
function given by
\begin{equation}
R^{H}(t,s):=\mathbb{E}\left[  B^{H}\left(  t\right)  B^{H}\left(  s\right)
\right]  =\frac{1}{2}(t^{2H}+s^{2H}-|t-s|^{2H}). \label{cov}%
\end{equation}
The parameter $H\ $characterizes all the important properties of the process.
In addition, to being self-similar with parameter $H$, which is evident from
the covariance function, fBm has correlated increments: in fact, from
(\ref{cov}) we get, as $n\rightarrow\infty$,
\begin{equation}
\mathbb{E}\left[  \left(  B^{H}\left(  n\right)  -B^{H}\left(  1\right)
\right)  B^{H}\left(  1\right)  \right]  =H\left(  2H-1\right)  n^{2H-2}%
+o\left(  n^{2H-2}\right)  ; \label{long}%
\end{equation}
when $H<1/2$, the increments are negatively correlated and the correlation
decays more slowly than quadratically; when $H>1/2$, the increments are
positively correlated and the correlation decays so slowly that they are not
summable, a situation which is commonly known as the long memory property. The
covariance structure (\ref{cov}) also implies%
\begin{equation}
\mathbb{E}\left[  \left(  B^{H}\left(  t\right)  -B^{H}\left(  s\right)
\right)  ^{2}\right]  =\left\vert t-s\right\vert ^{2H}; \label{canon}%
\end{equation}
this property shows that the increments of fBm are stationary and
self-similar; its immediate consequence for higher moments can be used, via
the so-called Kolmogorov continuity criterion, to imply that $B^{H}$ has paths
which are almost-surely ($H-\varepsilon$)-H\"{o}lder-continuous for any
$\varepsilon>0$.

It turns out that fBm is the \emph{only} continuous Gaussian process which is self-similar with
stationary increments. This constitutes an alternative definition of the process. However, there are other stochastic processes which, except for the
Gaussian character, share all the other properties above for $H>1/2$ (i.e. (\ref{cov}) which
implies (\ref{long}), the long-memory property, (\ref{canon}), and in many cases the
H\"{o}lder-continuity). In some models the Gaussian assumption may be implausible and in this
case one needs to use a different self-similar process with stationary increments to model the
phenomenon. Natural candidates are the Hermite processes: these non-Gaussian stochastic
processes appear as limits in the so-called Non-Central Limit Theorem (see \cite{BrMa},
\cite{DM}, \cite{Ta1}) and do indeed have all the properties listed above. While fBm can be
expressed as a Wiener integral with respect to the standard Wiener process, i.e. the integral of
a deterministic kernel w.r.t. a standard Brownian motion, the Hermite process of order $q\geq2$
is a $q$th iterated integral of a deterministic function with $q$ variables with respect to a
standard Brownian motion. When $q=2$, the Hermite process is called the Rosenblatt process. This
stochastic process typically appears as a limiting model in various applications such as unit
the root testing problem (see \cite{Wu}) or  semiparametric approach to hypothesis test (see
\cite{Hall}). On the other hand, since it
is non-Gaussian and self-similar with stationary increments, the Rosenblatt process can also be
an input in models where self-similarity is observed in empirical data which appears to be
non-Gaussian. The need of non-Gaussian self-similar processes in practice (for example in
hydrology) is mentioned in the paper \cite{Taqqu3} based on the study of stochastic modeling for
river-flow time series in \cite{LK}.  Recent interest in the Rosenblatt and other Hermite
processes, due in part to their non-Gaussian character, and in part for their independent
mathematical value, is evidenced by the following references: \cite{BN}, \cite{CNT},
 \cite{NNT}, \cite{T}, \cite{TV}.

 In this paper we will give a strong approximation result for the Rosenblatt process by means of transport processes. It is also interesting  from the theoretical point of view since all the approximation results for the Rosenblatt process known in the literature are in the weak sense (\cite{DM}, \cite{Ta1}).

 Our work is a natural extension of the strong approximation results for the Brownian motion and for the fractional Brownian motion. The study of the convergence of transport processes to the Brownian motion has a long history. We mention the works (\cite{CH}, \cite{GHR}, \cite{GoGr}) among others. More recently, due to the development of the stochastic analysis for fractional Brownian motion, the need of simulating the paths of this process led to the study of the strong approximation of the fBm by means. We refer to \cite{GaGoLe} for such an approximation in terms of transport processes and to \cite{E} or \cite{Sz} for related works.

 Our paper is organized as follows. In Section 2 we give some preliminares on multiple integrals and Malliavin derivatives. In section 3 we describe the approximating processes and prove the convergence to Rosenblatt process.

\section{Multiple Wiener-It\^o Integrals and Malliavin Derivatives}
 We start by introducing the elements from stochastic analysis that we will need in the paper. Consider ${\mathcal{H}}$ a real separable Hilbert space and $(B (\varphi), \varphi\in{\mathcal{H}})$ an isonormal Gaussian process on a probability space $(\Omega, {\cal{A}}, P)$, which is a centered Gaussian family of random variables such that $\mathbb{E}\left( B(\varphi) B(\psi) \right)  = \langle\varphi, \psi\rangle_{{\mathcal{H}}}$. Denote by $I_{n}$ the multiple stochastic integral with respect to
$B$ (see \cite{N}). This $I_{n}$ is actually an isometry between the Hilbert space ${\mathcal{H}}^{\odot n}$(symmetric tensor product) equipped with the scaled norm $\frac{1}{\sqrt{n!}}\Vert\cdot\Vert_{{\mathcal{H}}^{\otimes n}}$ and the Wiener chaos of order $n$ which is defined as the closed linear span of the random variables $H_{n}(B(\varphi))$ where $\varphi\in{\mathcal{H}}, \Vert\varphi\Vert_{{\mathcal{H}}}=1$ and $H_{n}$ is the Hermite polynomial of degree $n\geq 1$
\begin{equation*}
H_{n}(x)=\frac{(-1)^{n}}{n!} \exp \left( \frac{x^{2}}{2} \right)
\frac{d^{n}}{dx^{n}}\left( \exp \left( -\frac{x^{2}}{2}\right)
\right), \hskip0.5cm x\in \mathbb{R}.
\end{equation*}
The isometry of multiple integrals can be written as: for $m,n$ positive integers,
\begin{eqnarray}
\label{eqisometry}
\mathbb{E}\left(I_{n}(f) I_{m}(g) \right) &=& n! \langle f,g\rangle _{{\mathcal{H}}^{\otimes n}}\quad \mbox{if } m=n,\nonumber \\
\mathbb{E}\left(I_{n}(f) I_{m}(g) \right) &= & 0\quad \mbox{if } m\not=n.\label{iso}
\end{eqnarray}
It also holds that
\begin{equation*}
I_{n}(f) = I_{n}\big( \tilde{f}\big)
\end{equation*}
where $\tilde{f} $ denotes the symmetrization of $f$ defined by $\tilde{f}%
(x_{1}, \ldots , x_{n}) =\frac{1}{n!} \sum_{\sigma \in {\cal S}_{n}}
f(x_{\sigma (1) }, \ldots , x_{\sigma (n) } ) $.
\\\\
We recall the following hypercontractivity property  for the $L^{p}$ norm of a multiple
stochastic integral (see \cite[Theorem 4.1]{Major})
\begin{equation}
  \label{hyper}
  \bE \left| I_{m}(f) \right| ^{2m} \leq c_{m} \left( \bE I_{m}(f)^{2}
  \right)^{m}
\end{equation}
where $c_{m}$ is an explicit positive constant and $f\in
{\cal{H}}^{\otimes m}$.

\noindent In this paper we will use multiple stochastic integrals with respect to the Brownian motion (Bm) on $\mathbb{R}$ as introduced above. Note that the Brownian motion on the real line is an isonormal process and its underlying Hilbert space is ${\cal{H}}=L^{2}(\mathbb{R})$.
\\\\
\noindent For every $\frac{1}{2} \leq H <1$ the Rosenblatt process $(X^{H} _{t})_{t\in [0,T] }  $ could be defined as follows,
\begin{equation}\label{hermite}
X^{H}_{t} =c(H) I_{2} (g_{t}(\cdot ) )
\end{equation}
where for every $t\in [0,T]$
\begin{equation}\label{kernel}g_{t}(y_{1},y_{2})= \int_{y_{1}\vee y_{2}}^{t} (u-y_{1})_{+}^{\frac{H}{2}-1 }(u-y_{2})_{+}^{\frac{H}{2} -1 }du.
\end{equation}
The constant $d(H)$ is a normalizing constant which ensures that $\mathbb{E}(X^{H}_{t})^{2} =t^{2H}$ for every $t\in [0,T]$. This constant can be explicitly computed but it has no interest for our investigation. It can be proved that the process $X^{H}$ is self-similar with stationary increment and has the same covariance (\ref{cov}) as the fBm. Moreover it satisfies properties (\ref{long}) and (\ref{canon}).

\section{Strong convergence to the Rosenblatt process}
The Rosenblatt process $(X^{H} _{t}) _{t\in [0,T]} $ defined above can be written as an iterated double integral in the following way
\begin{equation}\label{rose}
X^{H}_{t}= c(H) \int _{\mathbb{R}} \int_{\mathbb{R}} \left( \int_{0} ^{t}   (s-x_{1} )_{+} ^{\frac{H}{2}-1}(s-x_{2} )_{+} ^{\frac{H}{2}-1} ds\right)dB(x_{1} ) dB(x_{2}), \hskip0.5cm t\in [0,T]
\end{equation}
where $B$ is a Wiener process on the whole real line and the Hurst parameter $H$ belongs to the interval $(\frac{1}{2} , 1)$. The process $X^{H}$ is $H$ self similar with stationary increments and it has the same covariance as the fractional Brownian motion.

We will separate $X^{H}$ into three terms. For every $t\in [0,T]$
\begin{eqnarray}
\label{defx1x4}
X^{H}_{t}&=& c(H) \biggl[\int_{-\infty} ^{0} \int _{-\infty} ^{0}  \left( \int_{0} ^{t}   (s-x_{1} ) ^{\frac{H}{2}-1}(s-x_{2} ) ^{\frac{H}{2}-1} ds\right)dB(x_{1} ) dB(x_{2})\biggr. \notag \\
&&+ \int_{-\infty} ^{0} \int_{0} ^{t}  \left( \int_{x_1} ^{t}   (s-x_{1} ) ^{\frac{H}{2}-1}(s-x_{2} )^{\frac{H}{2}-1} ds\right)dB(x_{1} ) dB(x_{2})\notag\\
&&+ \int_{0} ^{t} \int_{-\infty} ^{0}  \left( \int_{x_2} ^{t}   (s-x_{1} ) ^{\frac{H}{2}-1}(s-x_{2} ) ^{\frac{H}{2}-1} ds\right)dB(x_{1} ) dB(x_{2})\notag\\
&&\biggl.+ \int_{0} ^{t} \int_{0}^{t}  \left( \int_{x_{1} \vee x_{2}} ^{t}   (s-x_{1} ) ^{\frac{H}{2}-1}(s-x_{2} ) ^{\frac{H}{2}-1} ds\right)dB(x_{1} ) dB(x_{2})\biggr]\notag\\
&:=& X^{1,H}_{t}+ 2X^{2,H}_{t}+ X^{3, H}_{t},
\end{eqnarray}
note that the second and the third integrals are actually equal, for that reason the term $X^{2,H}$ appears twice . We will treat separately the third terms above since they have different behavior which comes from the singularity of the integral appearing in their expression.

\subsection{Transport processes}	

For each $n=1,2, \ldots$, let $(Z^{(n)}(t))_{ t\geq 0}$ be a process such that $Z^{(n)}(t)$ is the position on the real line at time $t$ of a particle moving as follows. It starts from $0$
with constant velocity $+n$ or $-n$, each with probability 1/2. It continues until a random time $\tau_1$ which is  exponentially distributed with parameter $n^2$, and at that time it switches from velocity $\pm n $ to $\mp n$ and continues for an additional independent random time $\tau_2-\tau_1$ which is again exponentially distributed with parameter $n^2$. At time $\tau_2$ it changes velocity as before, and so on.  This process is called a (uniform) transport process. Griego, Heath and Ruiz-Moncayo \cite{GHR} showed that $Z^{(n)}$ converges  to Brownian motion strongly and uniformly on bounded time intervals, and a rate of convergence was derived by Gorostiza and Griego in \cite{GoGr} as follows,

\begin{theorem}
\label{taproxmb}
There exist versions  of the transport processes $Z^{(n)}$ on the same probability space as a given Brownian motion  $(B_t)_{t\geq 0}$ such that for each $q>0$,
\begin{equation*}
	P\left(\sup_{a\leq t\leq b}|B_t- Z^{(n)}_t|> Cn^{-1/2}(\log n)^{5/2}\right)= o(n^{-q}) \  \ \  \text{as} \ n\to \infty,
\end{equation*}
 where $C$ is a positive constant depending on $a, b$ and $q$.
\end{theorem}

Let  $(X^{H} _{t}) _{t\in [0,T]} $ a Rosenblatt process. With $a<0$  fixed, we consider the following Bm's  constructed from the Bm $B$ in (\ref{rose}),
\begin{enumerate}

\item $\left( B_{1}(s) \right) _{s\in [0,T]} $ , the restriction of $B$ to the interval $[0,T]$.
    \item $\left(B_2(s)\right)_{a \leq s \leq 0}$, the restriction of $B$ to the interval $\left[a, 0\right]$.
    \item $B_3(s)=\begin{cases} sB(\frac{1}{s}) & \text{if}\  s\in \left[\frac{1}{a}, 0\right),\\
    0& \text{if}\  s=0. \end{cases}$
\end{enumerate}
Let us define now the transport processes that will intervene in our main results. By Theorem \ref{taproxmb}, there are three transport processes
\begin{equation}
\label{edeftranspro1}
	 (Z_1^{(n)}(s))_{0 \leq s \leq T}, \ \  (Z_2^{(n)}(s))_{a \leq s \leq 0}, \ \  \text{and} \ \ \ (Z_3^{(n)}(s))_{ \frac{1}{a}\leq s \leq 0},
\end{equation}
 such that  for each $q>0$,
\begin{equation}
\label{eq98}
    P\left(\sup_{b_i\leq t \leq c_i}|B_i(t) - Z_i^{(n)}(t)|> C^{(i)} n^{-1/2}(\log n )^{5/2}\right)= o(n^{-q}) \ \ \ \text{as}  \ \ n\to \infty,
\end{equation}
where $b_i, c_i$, $i=1,2,3$,  are the endpoints of the corresponding intervals, and $C^{(i)}$ is a positive constant depending on $b_i$, $c_i$ and $q$.

\subsection{Strong approximation}

We will approximate successively  each summand $X^{1,H}, X^{2,H}, X^{3,H}$ from (\ref{defx1x4}) in the strong sense by processes construct in terms of the transport processes $Z_{1}^{(n)}, Z^{(n)} _{2}, Z^{(n)}_{3}$ introduced above. Let us start with the summand $X^{1,H}$. Using Fubini theorem, we can express it as
\begin{eqnarray}
\label{eqdefx1}X^{1,H}_{t}&=&c(H) \int_{0}^{t} ds \int _{-\infty} ^{0} \int_{-\infty} ^{0} dB(x_{1}) dB(x_{2})  (s-x_{1} ) ^{\frac{H}{2}-1}(s-x_{2} ) ^{\frac{H}{2}-1}\notag\\
&=& c(H) \int_{0} ^{t} ds \left( \int_{-\infty} ^{0} (s-x) ^{\frac{H}{2}-1} dB(x) \right) ^{2}\\
&=& c(H) \int_{0} ^{t} ds \left(  Y^{1,H}_{s}\right) ^{2}, \hskip0.5cm t\in [0,T]\notag
\end{eqnarray}
where
\begin{equation}\label{y1H}
Y^{1,H}_{s}= \int_{-\infty} ^{0} (s-x) ^{\frac{H}{2}-1} dB(x), \hskip0.5cm s\in [0,T].
\end{equation}
\begin{remark}
Notice that integral $\int_{-\infty} ^{0} (s-x) ^{\frac{H}{2}-1} dB(x)$ is well-defined in $L^{2}(\Omega)$ as a Wiener integral for every $s>0$ since
\begin{eqnarray*}
\mathbb{E} \left( \int_{-\infty} ^{0} (s-x) ^{\frac{H}{2}-1} dB(x) \right) ^{2} &=& \int_{-\infty} ^{0}  (s-x) ^{H-2} dx = \frac{1}{1-H} s^{2H-1}.
\end{eqnarray*}
The situation will be different when we treat the summand $X^{3,H}$. This is one of the reasons to decompose the Rosenblatt process into several parts.
\end{remark}

Let $0<\max \left(  \frac{1-H/2}{3-2H}, \frac{2-H}{2H+2} \right)<\beta< 1/2$ be fixed (note that $\frac{1-H/2}{3-2H} <\frac{1}{2}$ since $H<1$ and $\frac{2-H}{2H+2}<\frac{1}{2}$ because $H>\frac{1}{2}$), denote in the sequel by
\begin{equation}
\label{eqdefepsilon}
\varepsilon_n=n^{-\frac{ \beta}{1-H/2}}
\end{equation}
and by
\begin{equation}
\label{eqdefalphan}
\alpha _{n}= n^{-(\frac{1}{2} -\beta ) }(\log n) ^{\frac{5}{2} }.
\end{equation}

We  will use the notation $$\left\|Y\right\|_{\infty, [a,b]}=\sup_{a\leq s \leq b}|Y_s|.$$
When the interval is of the form $[0,T]$ we will use the shorter notation $\left\|Y\right\|_{\infty, [0,T]}:= \left\|Y\right\|_{\infty, T}.$
\noindent We will denoted by $C$ a generic strictly positive constant that may depend on $a, T, H, p$ and may change from line to line.

\

Let us give a different expression for the process $Y^{1,H}$.

\begin{lemma}
\label{lemma1}
Let $Y^{1,H}$ be the process defined by (\ref{y1H}) and   $a<0$ fixed, then for every $s\in [0,T]$
\begin{eqnarray}
\label{eqdefY3}
Y^{1,H}_{s}&= & f_s(a)B_2(a)-\int_{1/a}^{-\varepsilon_n}\partial_xf_s\left(\frac{1}{u}\right)\frac{1}{u^3}B_3({u})du-\int_{-\varepsilon_n}^0\partial_xf_s\left(\frac{1}{u}\right)\frac{1}{u^3}B_3({u})du\notag \\
&&+\int_a^{-\varepsilon_n}f_s(x)dB_2(x) + \int_{-\varepsilon_n}^0[f_s(x)-f_s(x-\varepsilon_n)]dB_2(x)\notag\\
&&+ \int_{-\varepsilon_n}^0f_s(x-\varepsilon_n)dB_2(x)
\end{eqnarray}
where $\partial_xf_s$ denotes the derivative of the function $f_{s}(x)= (s-x) ^{H/2-1}, s>x$ with respect to its second variable (even when this second variable is not denoted by $x$).
\end{lemma}
{\bf Proof: } We can write, for every $t\in [0,T]$
\begin{equation}
\label{eqdefY1}
    Y^{1,H}_{s}= \int_{-\infty}^{a}f_s(x)dB(x)+ \int_{a}^{0}f_s(x)dB(x).
\end{equation}
 We express the first Wiener integral above as an integral with respect to $ds$.  Since  by the  H\"older continuity of $B$,
$$\lim_{b\to -\infty}f_s(b)B(b)=0,$$
by integration by parts and putting $x=1/u$,
\begin{eqnarray}
\label{eqdefY2}
\int_{-\infty}^{a}f_s(x)dB(x)&=& f_s(a)B(a)-\int_{-\infty}^{a}\partial_xf_s(x)B(x)dx\nonumber \\
&= &f_s(a)B(a)-\int_{1/a}^{0}\partial_xf_s\left(\frac{1}{u}\right)\frac{1}{u^2}B\left(\frac{1}{u}\right)du\nonumber\\
&=& f_s(a)B_2(a)-\int_{1/a}^{0}\partial_xf_s\left(\frac{1}{u}\right)\frac{1}{u^3}B_3({u})du.
\end{eqnarray}
By (\ref{eqdefY1}) and (\ref{eqdefY2}), for every $s\in [0,T]$ we have the result.
\qed

\vskip0.2cm

We first approximate the process $(Y^{1,H}_{s})_{s\in [0,T]} $ (in the strong sense (\ref{eq98})) by  stochastic processes constructed from transport processes. Basically, in the expression of $Y^{1,H}$, we replace the Brownian motions by their corresponding transport processes.
The approximating processes to $Y^{1,H}$ is defined as

\begin{eqnarray}
\label{eqdefYn}
Y^{1,H, n}_{s}&=& f_s(a)Z_2^{(n)}(a)-\int_{1/a}^{-\varepsilon_n}\partial_xf_s\left(\frac{1}{u}\right)\frac{1}{u^3}Z_3^{(n)}({u})du\notag +\int_a^{-\varepsilon_n}f_s(x)dZ_2^{(n)}(x)\notag\\
&& + \int_{-\varepsilon_n}^0f_s(x-\varepsilon_n)dZ_2^{(n)}(x), \hskip0.5cm s\in [0,T].
\end{eqnarray}

We state the result concerning the approximation of $Y^{1,H}$. Its proof follows the ideas of the proofs in \cite{GaGoLe} but the context is technically more complex. Note that the singularity of the integrand $(s-x)^{\frac{H}{2} -1}$ at $s=x$ does not allows to use directly the results in \cite{GaGoLe} and the arguments of the proofs must be adapted to fit in our context.
\begin{prop}
\label{prop1}
Let $Y^{1,H}$ and $Y^{1,H,n}$ be the processes defined by (\ref{y1H}) and (\ref{eqdefYn}), respectively and let $\alpha_n$ given by (\ref{eqdefalphan}). Then for each $q>0$ and each $\beta$ such that $0<\frac{1-H/2}{3-2H}< \beta < \frac{1}{2}$,
\begin{equation}
        P\left(\sup_{0\leq s \leq T}s^{1-H/2}\left|Y^{1,H}_{s} - Y^{1,H, n}_{s}\right|> C\alpha_n\right)=o(n^{-q}) \ \ \text{as} \ n\to \infty.
\end{equation}
\end{prop}


{\bf Proof: } From (\ref{eqdefY3}) and  Lemma \ref{lemma1} we have
\begin{eqnarray*}
    &&|Y^{1,H}_{t}-Y^{1,H,n}_{t}|\leq \Biggl\{\Biggl.\left|f_t(a)B_2(a)-f_t(a)Z_2^{(n)}(a)\right|\\
    &&+\left|\int_{1/a}^{-\varepsilon_n}\partial_xf_s\left(\frac{1}{u}\right)\frac{1}{u^3}B_3({u})du-\int_{1/a}^{-\varepsilon_n}\partial_xf_s\left(\frac{1}{u}\right)\frac{1}{u^3}Z_3^{(n)}(u)du\right|\\
    &&+\left|\int_{-\varepsilon_n}^0\partial_xf_s\left(\frac{1}{u}\right)\frac{1}{u^3}B_3({u})du\right|+\left|\int_a^{-\varepsilon_n}f_s(x)dB_2(x)-\int_a^{-\varepsilon_n}f_s(x)dZ_2^{(n)}(x)\right|\\
&&+\left|\int_{-\varepsilon_n}^0f_s(x-\varepsilon_n)dB_2(x)-\int_{-\varepsilon_n}^0f_s(x-\varepsilon_n)dZ_2^{(n)}(x)\right|+ \left|\int_{-\varepsilon_n}^0[f_s(x)-f_s(x-\varepsilon_n)]dB_2(x)\right|\Biggl.\Biggr\}.
\end{eqnarray*}
By  Lemmas \ref{lemB1}, \ref{lemA1}, \ref{lemD1}, \ref{lemE1}, \ref{lemF1} and \ref{lemG1} below we have the result.
\qed

\begin{lemma}
\label{lemB1}Let $Z_2^{(n)}$ be the process defined by (\ref{edeftranspro1}). Then for each $q>0$ there is $C>0$ such that
\begin{equation}
\label{eqB1}
I_1: =P\left(\sup_{0\leq s \leq T}\left| f_s(a)B_2(a) -f_s(a)Z_2^{(n)}(a)\right|> C\alpha_n \right)= o(n^{-q})\ \ \ \text{as} \ \   n\to \infty.
\end{equation}
\end{lemma}
{\bf Proof: } It holds, for fixed $a<0$,
\begin{eqnarray*}
\left| f_s(a)B_2(a) - f_s(a)Z_2^{(n)}(a))\right|&\leq & \|B_2 - Z_2^{(n)}\|_{\infty,[a,0]}(s-a)^{H/2-1}\\
&\leq&  \|B_2 - Z_2^{(n)}\|_{\infty,[a,0]}(-a)^{H/2-1}
\end{eqnarray*}
then (recall that $C$ is a generic strictly positive constant that may depend on $a, T, H$)
\begin{eqnarray*}
I_1&\leq &P\left(\|B_2 - Z_2^{(n)}\|_{\infty,[a,0]}(-a)^{H/2-1}> C\alpha_n \right)\\
&\leq &P\left(\|B_2 - Z_2^{(n)}\|_{\infty,[a,0]}> C\alpha_n \right)=o(n^{-q}).
\end{eqnarray*}\qed

\begin{remark}
The conclusion of Lemma \ref{lemB1}  is clearly true if we add  the factor $s^{1-\frac{H}{2} } $ after the supremum. This remark is also available for the following lemmas and we will not mention it at each time.
\end{remark}

\begin{lemma}
\label{lemA1} Let $Z_3^{(n)}$ be the process defined by (\ref{edeftranspro1}). Then for each $q>0$,
\begin{eqnarray}
\label{eqA1}
I_2&:=&P\left(\sup_{0\leq s \leq T}s^{1-H/2}\left|\int_{1/a}^{-\varepsilon_n}\partial_xf_s\left(\frac{1}{u}\right)\frac{1}{u^3}B_3({u})du-\int_{1/a}^{-\varepsilon_n}\partial_xf_s\left(\frac{1}{u}\right)\frac{1}{u^3}Z_3^{(n)}({u})du\right|> C\alpha_n  \right)\notag\\
&= &o(n^{-q}) \ \ \ \text{as} \ \   n\to \infty.
\end{eqnarray}
\end{lemma}
{\bf Proof: }Putting $z=1/u$ and $w=s-z$,
\begin{eqnarray*}
&&\left|\int_{1/a}^{-\varepsilon_n}\partial_xf_s\left(\frac{1}{u}\right)\frac{1}{u^3}B_3({u})du-\int_{1/a}^{-\varepsilon_n}\partial_xf_s\left(\frac{1}{u}\right)\frac{1}{u^3}Z_3^{(n)}({u})du\right|\\
&&\leq \|B_3 - Z_3^{(n)}\|_{\infty,[1/a,0]}\int_{1/a}^{-\varepsilon_n}\left|\partial_xf_s\left(\frac{1}{u}\right)\frac{1}{u^3}B_3({u})\right|du\\
&&=\|B_3 - Z_3^{(n)}\|_{\infty,[1/a,0]}(1-H/2)\int_{1/a}^{-\varepsilon_n}\frac{1}{(-u)^3}(s-1/u)^{H/2-2}du\\
&&=\|B_3 - Z_3^{(n)}\|_{\infty,[1/a,0]}(1-H/2)\int_{-1/\varepsilon_n}^a(-z)(s-z)^{H/2-2}dz\\
&&=\|B_3 - Z_3^{(n)}\|_{\infty,[1/a,0]}(1-H/2)\int^{s+1/\varepsilon_n}_{s-a}(w-s)(w)^{H/2-2}dw\\
&&\leq\|B_3 - Z_3^{(n)}\|_{\infty,[1/a,0]}(1-H/2)\int^{s+1/\varepsilon_n}_{s-a}w^{H/2-1}dw\\
&&=\|B_3 - Z_3^{(n)}\|_{\infty,[1/a,0]}\frac{1-H/2}{H/2}[(s+1/\varepsilon_n)^{H/2}-(s-a)^{H/2}]\\
&&\leq\|B_3 - Z_3^{(n)}\|_{\infty,[1/a,0]}\frac{1-H/2}{H/2}(T+1/\varepsilon_n)^{H/2}\\
&&\leq\|B_3 - Z_3^{(n)}\|_{\infty,[1/a,0]}\frac{1-H/2}{H/2}2^{H/2}(T^{H/2}+(\varepsilon_n)^{-H/2}),
\end{eqnarray*}
then
\begin{eqnarray*}
I_2&&\leq P\left(\|B_3 - Z_3^{(n)}\|_{\infty,[1/a,0]}\frac{1-H/2}{H/2}2^{H/2}(T^{H/2}+(\varepsilon_n)^{-H/2})>C\alpha_n\right)\notag\\
&&\leq P\left(\|B_3 - Z_3^{(n)}\|_{\infty,[1/a,0]}\frac{1-H/2}{H/2}2^{H/2}T^{H/2}>C\alpha_n\right)\\
&&+P\left(\|B_3 - Z_3^{(n)}\|_{\infty,[1/a,0]}\frac{1-H/2}{H/2}2^{H/2}T^{1-H/2}(\varepsilon_n)^{-H/2}>C\alpha_n\right)\notag\\
&&\leq P\left(\|B_3 - Z_3^{(n)}\|_{\infty,[1/a,0]}>C\alpha_n\right)\\
&&+P\left(\|B_3 - Z_3^{(n)}\|_{\infty,[1/a,0]}>C(\varepsilon_n)^{H/2}\alpha_n\right)\notag\\
&&\leq o(n^{-q})+ P\left(\|B_3 - Z_3^{(n)}\|_{\infty,[1/a,0]}>Cn^{-1/2+\beta(1-H)/(1-H/2)}(\log n)^{5/2}\right)\notag \\
&&=o(n^{-q}).
\end{eqnarray*}\qed

\vskip0.2cm

The following lemma explains one of the  conditions imposed on $\beta $ in the statement of Proposition \ref{prop1}. Another restriction comes from Proposition 4 later.
\begin{lemma}
\label{lemD1} Let $(1-H/2)/(3-2H)< \beta < \frac{1}{2}$. Then for each $q>0$
\begin{equation}
\label{eqD1}
I_3:=P\left(\sup_{0\leq s \leq T}\left|\int_{-\varepsilon_n}^0\partial_xf_s\left(\frac{1}{u}\right)\frac{1}{u^3}B_3({u})du\right|> \alpha_n\right)=o(n^{-q})\ \ \ \text{as} \ \  n\to \infty.
\end{equation}
\end{lemma}
{\bf Proof: }For every $s\in [0,T]$ and $u\in [-\eps _{n}, 0)$  we can write
\begin{eqnarray}
\label{eqD2}\left|\frac{1}{u^3}\partial_xf_s\left(\frac{1}{u}\right)\right|&&=\left|\frac{1}{u^3}(H/2-1).-(s-1/u)^{H/2-2}\right|
=\frac{1-H/2}{(-u)^3}\cdot\left(\frac{1-us}{-u}\right)^{H/2-2}\notag\\
&& \leq (1-H/2)(-u)^{-1-H/2}
\end{eqnarray}

By (\ref{eqD2}) and the pathwise H\"older continuity of the Bm $B_3$ there exists a random variable $Y$ (having all its moments finite) such that for any  $\gamma<1/2-H/2$,
\begin{eqnarray*}
\left|\int_{-\varepsilon_n}^0\partial_xf_s\left(\frac{1}{u}\right)\frac{1}{u^3}B_3({u})du\right|&&\leq Y\int_{-\varepsilon_n}^0 (1-H/2)(-u)^{-1-H/2}(-u)^{1/2-\gamma}du \\
&&=Y\frac{1-H/2}{1/2-H/2-\gamma}(\varepsilon_n)^{1/2-H/2-\gamma}.
\end{eqnarray*}
By Chebyshev's inequality, for $r>0$,
\begin{eqnarray*}
I_3&&\leq P\left(Y\frac{1-H/2}{1/2-H/2-\gamma}n^{-\beta\frac{1/2-H/2-\gamma}{1-H/2}}>\alpha_n\right)\\
&&=P\left(CY>n^{\kappa}(\log n)^{5/2}\right)\leq \frac{\mathbb{E}(|\tilde{C}Y|^r)}{n^{r\kappa}(\log n)^{r5/2}},
\end{eqnarray*}
where $\kappa=-(1/2-\beta)+ \beta(1/2-H/2-\gamma)/(1-H/2)$. Taking
$$(1-H/2)/(3-2H)<(1-H/2)/(3-2H-\gamma)<\beta<1/2,$$
 then $\kappa>0$. For $q>0$ there is $r>0$ such that $q<r\kappa$, then
\begin{equation*}
	\lim_{n\to \infty}n^qI_3=0.
\end{equation*}
\qed

\begin{lemma}
\label{lemE1} Let $Z_2^{(n)}$ be the process defined by (\ref{edeftranspro1}). Then for each $q>0$,
\begin{eqnarray}
\label{eqE1}
I_4&:=&P\left(\sup_{0\leq s \leq T}\left|\int_a^{-\varepsilon_n}f_s(x)dB_2(x)-\int_a^{-\varepsilon_n}f_s(x)dZ_2^{(n)}(x)\right|
> C\alpha_n \right)\notag\\
&=& o(n^{-q}) \ \ \ \text{as} \ n\to \infty
\end{eqnarray}
\end{lemma}
{\bf Proof: }By integration by parts,
\begin{eqnarray*}
\int_a^{-\varepsilon_n}f_s(x)dB_2(x)= f_s(-\varepsilon_n)B_2(-\varepsilon_n)-f_s(a)B_2(a)-\int_a^{-\varepsilon_n}(1-H/2)(s-x)^{H/2-2}B_2(x)dx
\end{eqnarray*}
and
\begin{eqnarray*}
\int_a^{-\varepsilon_n}f_s(x)dZ_2^{(n)}(x)&=& f_s(-\varepsilon_n)Z_2^{(n)}(-\varepsilon_n)-f_s(a)Z_2^{(n)}(a)\\
&&-\int_a^{-\varepsilon_n}(1-H/2)(s-x)^{H/2-2}Z_2^{(n)}(x)dx
\end{eqnarray*}
then
\begin{eqnarray*}
&&\left|\int_a^{-\varepsilon_n}f_s(x)dB_2(x)-\int_a^{-\varepsilon_n}f_s(x)dZ_2^{(n)}(x)\right|\\
&&\leq \|B_2 - Z_2^{(n)}\|_{\infty,[a,0]}\left\{f_s(-\varepsilon_n)+f_s(a)+\int_a^{-\varepsilon_n}(1-H/2)(s-x)^{H/2-2}dx\right\}\\
&&\leq \|B_2 - Z_2^{(n)}\|_{\infty,[a,0]}\left\{(s+\varepsilon_n)^{H/2-1}+(s-a)^{H/2-1}+(s+\varepsilon_n)^{H/2-1}-(s-a)^{H/2-1}\right\}\\
&&= \|B_2 - Z_2^{(n)}\|_{\infty,[a,0]}2(s+\varepsilon_n)^{H/2-1}\\
&&\leq \|B_2 - Z_2^{(n)}\|_{\infty,[a,0]}2(\varepsilon_n)^{H/2-1}= \|B_2 - Z_2^{(n)}\|_{\infty,[a,0]}2n^{\beta}
\end{eqnarray*}
Consequently,
\begin{eqnarray*}
I_4&&\leq P\left(\|B_2 - Z_2^{(n)}\|_{\infty,[a,0]}2n^{\beta}> C\alpha_n \right)\\
&&=P\left(\|B_2 - Z_2^{(n)}\|_{\infty,[a,0]}> Cn^{-1/2}(\log n)^{5/2} \right)=o(n^{-q})
\end{eqnarray*}
\qed

\begin{lemma}
\label{lemF1} Let $Z_2^{(n)}$ be the process defined by (\ref{edeftranspro1}). Then for each $q>0$,
\begin{eqnarray}
\label{eqF1}
I_5&&:=P\left(\sup_{0\leq s \leq T}\left|\int_{-\varepsilon_n}^0f_s(x-\varepsilon_n)dB_2(x)-\int_{-\varepsilon_n}^0f_s(x-\varepsilon_n)dZ_2^{(n)}(x)\right|
> C\alpha_n \right)\notag\\
&&= o(n^{-q}) \ \ \ \text{as} \ \   n\to \infty.
\end{eqnarray}
\end{lemma}
{\bf Proof: } By integration by parts as before and taking into account that $B_{2}(0)= Z^{(n)}_{2}(0)=0$ we can express the two integrals in the statement as
\begin{eqnarray*}
\int_{-\varepsilon_n}^0f_s(x-\varepsilon_n)dB_2(x)= -f_s(-2\varepsilon_n)B_2(-\varepsilon_n)-\int_{-\varepsilon_n}^0(1-H/2)(s+\varepsilon_n-x)^{H/2-2}B_2(x)dx
\end{eqnarray*}
and
\begin{eqnarray*}
\int_{-\varepsilon_n}^0f_s(x-\varepsilon_n)dZ_2^{(n)}(x)&=& -f_s(-2\varepsilon_n)Z_2^{(n)}(-\varepsilon_n)\\
 &&-\int_{-\varepsilon_n}^0(1-H/2)(s+\varepsilon_n-x)^{H/2-2}Z_2^{(n)}(x)dx.
\end{eqnarray*}
then
\begin{eqnarray*}
&&\left|\int_{-\varepsilon_n}^0f_s(x-\varepsilon_n)dB_2(x)-\int_{-\varepsilon_n}^0f_s(x-\varepsilon_n)dZ_2^{(n)}(x)\right|\\
&&\leq \|B_2 - Z_2^{(n)}\|_{\infty,[a,0]}\left\{f_s(-2\varepsilon_n)+\int_{-\varepsilon_n}^0(1-H/2)(s+\varepsilon_n-x)^{H/2-2}dx\right\}\\
&&\leq \|B_2 - Z_2^{(n)}\|_{\infty,[a,0]}\left\{(s+2\varepsilon_n)^{H/2-1}+(s+\varepsilon_n)^{H/2-1}-(s+2\varepsilon_n)^{H/2-1}\right\}\\
&&= \|B_2 - Z_2^{(n)}\|_{\infty,[a,0]}(s+\varepsilon_n)^{H/2-1}\leq \|B_2 - Z_2^{(n)}\|_{\infty,[a,0]}(\varepsilon_n)^{H/2-1}\\
&&= \|B_2 - Z_2^{(n)}\|_{\infty,[a,0]}n^{\beta}
\end{eqnarray*}
Consequently,
\begin{eqnarray*}
I_5&&\leq P\left(\|B_2 - Z_2^{(n)}\|_{\infty,[a,0]}n^{\beta}> C\alpha_n \right)\\
&&=P\left(\|B_2 - Z_2^{(n)}\|_{\infty,[a,0]}> Cn^{-1/2}(\log n)^{5/2} \right)=o(n^{-q})
\end{eqnarray*}
\qed

Finally, we prove our last auxiliary approximation result. Here we need to add the factor $s^{1-\frac{H}{2}}$ which appears in Proposition \ref{prop1}. This is due to the singularity of the derivative of $f_{s}(x)$ with respect to $x$.
$s^{1-\frac{H}{2}}$.
\begin{lemma}
\label{lemG1}  Let $(1-H/2)/(3-2H)< \beta < \frac{1}{2}$. Then for each $q>0$
\begin{eqnarray}
\label{eqG1}
I_6&:=&P\left(\sup_{0\leq s \leq T}s^{1-\frac{H}{2}}\left|\int_{-\varepsilon_n}^0[f_s(x)-f_s(x-\varepsilon_n)]dB_2(x)\right|
> C\alpha_n \right)\notag\\
&=& o(n^{-q}) \ \ \ \text{as} \ \   n\to \infty.
\end{eqnarray}
\end{lemma}
{\bf Proof: }We first write the difference $f_s(x)-f_s(x-\varepsilon_n)$ as an integral and we use Fubini theorem. We obtain
\begin{eqnarray*}
&&\int_{-\varepsilon_n}^0[f_s(x)-f_s(x-\varepsilon_n)]dB_2(x)=\int_{-\varepsilon_n}^0\int_{s-x}^{s+\varepsilon_n-x}(1-H/2)u^{H/2-2}dudB_2(x)\\
&&=(1-H/2)\left[\int^{s+\varepsilon_n}_su^{H/2-2}\int_{s-u}^{0}dB_2(x)du+\int^{s+2\varepsilon_n}_{s+\varepsilon_n}u^{H/2-2}\int_{-\varepsilon_n}^{s+\varepsilon_n-u}dB_2(x)du\right]\\
&&=(1-H/2)\left[\int^{s+\varepsilon_n}_su^{H/2-2}[B_2(0)-B_2(s-u)]du\right.\\
&&\left.+\int^{s+2\varepsilon_n}_{s+\varepsilon_n}u^{H/2-2}[B_2(s+\varepsilon_n-u)-B_2(-\varepsilon_n)]du\right]
\end{eqnarray*}
The H\"older continuity of the Wiener process  $B_2$ implies for every $0<\gamma <\frac{1}{2}$
\begin{eqnarray*}
&&\left|\int_{-\varepsilon_n}^0[f_s(x)-f_s(x-\varepsilon_n)]dB_2(x)\right|\\
&&\leq (1-H/2)Y\left[\int^{s+\varepsilon_n}_su^{H/2-2}[u-s]^{1/2-\gamma}du
+\int^{s+2\varepsilon_n}_{s+\varepsilon_n}u^{H/2-2}[s+2\varepsilon_n-u]^{1/2-\gamma}du\right]\\
&&\leq (1-H/2)Y\left[\int^{s+\varepsilon_n}_su^{H/2-2}[\varepsilon_n]^{1/2-\gamma}du
+\int^{s+2\varepsilon_n}_{s+\varepsilon_n}u^{H/2-2}[\varepsilon_n]^{1/2-\gamma}du\right]\\
&&\leq (1-H/2)Y\left[\int_s^{s+2\varepsilon_n}u^{H/2-2}\varepsilon_n^{1/2-\gamma}du\right]\\
&&\leq Y\varepsilon_n^{1/2-\gamma}s^{H/2-1}
\end{eqnarray*}
and consequently
\begin{eqnarray*}
P\left(\sup_{0\leq s \leq T}s^{1-\frac{H}{2}}\left|\int_{-\varepsilon_n}^0[f_s(x)-f_s(x-\varepsilon_n)]dB_2(x)\right|
> C\alpha_n \right) \leq P\left(CY\varepsilon_n^{1/2-\gamma}>\alpha_n \right)
\end{eqnarray*}
and the result follows by analogous arguments as in proof of Lemma \ref{lemD1}.
\qed

\begin{remark}
\label{remark2}
In particular Proposition \ref{prop1} implies that
\begin{equation*}
P\left( \limsup_{n} \Big\{ \sup_{0\leq s \leq T}s^{1-\frac{H}{2}}\left| Y^{1,H}_s -Y^{1,H, n}_s \right|> C\alpha_n \Big\} \right) =0
\end{equation*}
by using Borel-Cantelli lemma.
\end{remark}
\vskip0.3cm

We finish the strong  approximation of the term $X^{1,H}$ appearing in the decomposition of the Rosenblatt process $X^{H}$ in (\ref{defx1x4}). By (\ref{eqdefx1}), we define for every $n$ and $Y^{1,H, n }$ given by (\ref{eqdefYn})
\begin{equation}
\label{x1n}
X^{1,H,n} = c(H) \int_{0} ^{t} (Y^{1,H, n }_{s} ) ^{2} ds.
\end{equation}

We have the following.

\begin{prop}\label{p2}
Let $X^{1,H} $ be given by (\ref{eqdefx1}) and $\beta\in \left(\frac{1-H/2}{3-2H}, \frac{1}{2}\right) $ fixed. Define $X^{1,H,n} $ by (\ref{x1n}). Then for any $\gamma$ such that $0<\gamma<\beta$ and $\beta+\gamma<1/2$,
\begin{equation*}
P \left( \limsup_{n\to \infty} \{ \Vert X^{1,H,n}-X^{1,H} \Vert _{\infty , T} \geq C n^{-(1/2-\beta-\gamma)}(\log n)^{5/2}\} \right) =0.
\end{equation*}
\end{prop}
{\bf Proof: } Using the fact that $A^{2}-B^{2}= (A-B)^{2}+ 2B(A-B)$ we can write, for every $t\in [0,T]$
\begin{eqnarray*}
X^{1,H,n}_{t}-X^{1,H}_{t}&=& c(H)\int_{0}^{t} ((Y^{1,H, n} _{s}) ^{2} -(Y^{1,H}_{s}) ^{2}) ds \\
&=&c(H)\int_{0} ^{t} \left[(Y^{1,H, n} _{s}-Y^{1,H}_{s})^{2} + 2Y^{1,H}_{s}(Y^{1,H,n}_{s}-Y^{1,H}_{s}) \right]ds
\end{eqnarray*}
and hence
\begin{eqnarray*}
\sup _{t\in [0,T] }|X^{1,H,n}_{t}-X^{1,H}_{t}|
&\leq & c(H)\sup_{ t\in [0,T]} \int_{0} ^{t}  \left|(Y^{1,H,n} _{s}-Y^{1, H}_{s})^{2} + 2Y^{1,H}_{s}(Y^{1,H,n}_{s}-Y^{1,H}_{s}) \right| ds\\
&=& c(H)\int_{0} ^{T} \left|(Y^{1,H,n} _{s}-Y^{1,H}_{s})^{2} + 2Y^{1,H}_{s}(Y^{1,H,n}_{s}-Y^{1,H}_{s}) \right| ds\\
&\leq & c(H)\int_{0}^{T} (Y^{1,H,n} _{s}-Y^{1,H}_{s})^{2}ds + 2c(H) \int_{0}^{T} |Y^{1,H}_{s}| \left| Y^{1,H,n}_{s}-Y^{1,H}_{s}\right| ds \\
&\leq & C \sup_{s\in [0,T] } s^{2-H}(Y^{1,H,n} _{s}-Y^{1,H}_{s})^{2}\\
&+& 2C\int_{0}^{T}s^{\frac{H}{2}-1} \left|Y^{1,H}_{s}\right| ds \sup_{s\in [0,T] }s^{1-\frac{H}{2}} |Y^{1,H,n} _{s}-Y^{1,H}_{s}|.
\end{eqnarray*}
We used above the trivial inequality $P(\vert X \vert ^{2} \geq C\alpha _{n}) \leq P(\vert X \vert  \geq C\alpha _{n})$ for any random variable $X$.  We will get ($C$ denoted a generic strictly positive constant depending on $T,H$ that may change from line to line) by Proposition \ref{prop1}
\begin{eqnarray}
\label{eqaproxx1}
&&P \left( \Vert X^{1,H, n}-X^{1,H}\Vert _{\infty , T} > C  n^{-(1/2-\beta-\gamma)}(\log n)^{5/2} \right)\notag\\
&\leq & P \left(  \sup_{ s\in [0,T]}s^{2-H}| Y^{1,H,n}_{s}-Y^{1,H}_{s} | ^{2}  > C n^{-(1/2-\beta-\gamma)}(\log n)^{5/2} \right) \notag\\
&+& P  \left( \sup_{ s\in [0,T]} s^{1-\frac{H}{2}}| Y^{1,H,n}_{s}-Y^{1,H}_{s} |\int_{0}^{T} s^{\frac{H}{2}-1}|Y^{1,H}_{s}|ds    > C n^{-(1/2-\beta-\gamma)}(\log n)^{5/2} \right) \notag\\
&=&o(n^{-q}) +  P  \left( \int_{0}^{T} s^{\frac{H}{2}-1}|Y^{1,H}_{s}| ds \sup_{ s\in [0,T]} s^{1-\frac{H}{2}}| Y^{1,H,n}_{s}-Y^{1,H}_{s} |  > C n^{-(1/2-\beta-\gamma)}(\log n)^{5/2} \right) .\notag\\
\end{eqnarray}

We apply Lemma \ref{lims} below  with $$A=\int_{0} ^{T}s^{\frac{H}{2}-1} |Y^{1,H}_{s}| ds \ \ \ \text{and} \ \ \ \Gamma=\sup_{0\leq s \leq T}s^{1-\frac{H}{2}}\left| Y^{1,H}_s -Y^{1,H, n}_s \right|.$$ We note first that
\begin{eqnarray*}
\mathbb{E} \int_{0} ^{T}s^{\frac{H}{2}-1} |Y^{1,H}_{s}| ds&\leq & c(H)\int_{0}^{T} ds  \ s^{\frac{H}{2}-1}\left( \mathbb{E} (Y^{1,H}_{s} )^{2} \right) ^{\frac{1}{2}}ds \\
&\leq & c(H) \int_{0}^{T} ds \ s^{\frac{H}{2}-1}\left( \int_{-\infty}^{0} (s-x) ^{H-2} dx \right) ^{\frac{1}{2}}\\
&=& c(H) \int_{0}^{T}  ds \ s^{\frac{H}{2}-1} s^{\frac{H-1}{2}} =c(H) T^{H-\frac{1}{2}}
\end{eqnarray*}
and thus the random variable $A$ is almost surely finite. We obtain by (\ref{eqaproxx1}) and Remark \ref{remark2}
\begin{eqnarray*}
&&P \left(\limsup \{ \Vert X^{1,H,n}-X^{1,H} \Vert _{\infty ,T} > Cn^{-(1/2-\beta- \gamma)}(\log n)^{5/2} \} \right) \\
&\leq& P \left( \limsup_{n\to \infty}  \{ A \sup_{0\leq s \leq T}s^{1-\frac{H}{2}}\left| Y^{1,H}_s -Y^{1,H, n}_s \right|> C n^{-(1/2-\beta- \gamma)}(\log n)^{5/2} \}\right)\\
&\leq &P\left(\limsup_{n\to \infty}  \{  \sup_{0\leq s \leq T}s^{1-\frac{H}{2}}\left| Y^{1,H}_s -Y^{1,H, n}_s \right|> C \alpha _{n} \} \right) =0.
\end{eqnarray*}
\qed

The following lemma has been used in the proof of Proposition \ref{p2}.

\begin{lemma}\label{lims}
Let $A$ and $\Gamma$ be random variables with $A$ an almost surely finite. Then for every $\gamma >0$
\begin{eqnarray*}
&&P\left( \limsup_{n\to \infty}  \{ A\Gamma > C n^{-(1/2-\beta-\gamma)}(\log n)^{5/2}\} \right)\\
&\leq& P\left(\limsup_{n\to \infty}  \{  \Gamma>C n^{-(1/2-\beta)}(\log n)^{5/2}\} \right)
\end{eqnarray*}
with $C$ a generic strictly positive constant.
\end{lemma}
{\bf Proof: }
We prove the following inclusion
\begin{eqnarray*}
&&\limsup_{n\to \infty}\left\{A\Gamma> n^{-(1/2-\beta-\gamma)}(\log n)^{5/2}\right\}\subseteq \limsup_{n\to \infty}\left\{\Gamma>n^{-(1/2-\beta)}(\log n)^{5/2}\right\}.
\end{eqnarray*}
Since
\begin{eqnarray*}
\omega \in& &	\limsup_{n\to \infty}\left\{A\Gamma>n^{-(1/2-\beta- \gamma)}(\log n)^{5/2}\right\}=\limsup_{n\to \infty}\left\{(An^{-\gamma})\Gamma>n^{-(1/2-\beta)}(\log n)^{5/2}\right\}\\
&=&\cap_{n=1}^{\infty}\cup_{k=n}^{\infty}\left\{(Ak^{-\gamma})\Gamma>k^{-(1/2-\beta)}(\log k)^{5/2}\right\},
\end{eqnarray*}
then for all $n\geq 1$,
$$\omega \in \cup_{k=n}^{\infty}\left\{(Ak^{-\gamma})\Gamma>k^{-(1/2-\beta)}(\log k)^{5/2}\right\},$$ and since $A$ is an almost surely finite random variable, there is $\hat{N}=\hat{N}(\omega)$ such that for all $n\geq \hat{N}$, $An^{-\gamma}<1$, then $\omega \in \cup_{k=n}^{\infty}\left\{\Gamma>k^{-(1/2-\beta)}(\log k)^{5/2}\right\}$
and the conclusion follows easily.
\qed

\vskip0.3cm

 \vskip0.3cm

Let us handle now the term $X^{3,H}$ appearing in (\ref{defx1x4}). We will decompose it as follows:
\begin{eqnarray}
\label{eqdescx4}
X^{3,H}_{t}&=& c(H) \int_{0}^{t}\int_{0}^{t} dB(x_{1}) dB(x_{2}) \int _{(x_{1}\vee x_{2} +\varepsilon _{n})\wedge t } ^{t}ds (s-x_{1} )^{\frac{H}{2}-1}(s-x_{2} )^{\frac{H}{2}-1}\notag\\
&&+ c(H) \int_{0}^{t}\int_{0}^{t} dB(x_{1}) dB(x_{2}) \int _{x_{1}\vee x_{2}}^{(x_{1}\vee x_{2} +\varepsilon _{n})\wedge t } ds (s-x_{1} )^{\frac{H}{2}-1}(s-x_{2} )^{\frac{H}{2}-1}\notag\\
&:=& c(H)(F^{n}_{t}+ G_{t}^{n})
\end{eqnarray}
where $\varepsilon _{n}$ is given by (\ref{eqdefepsilon}).
\begin{remark}
For the term $X^{3,H}$ we cannot use Fubini theorem (as in the case of $X^{1,H}$) because $\int_{0} ^{s} (s-x) ^{\frac{H}{2}-1} dB(x)$ is not defined as a Wiener integral in $L^{2}(\Omega)$ since the function $(s-x) ^{H-2}$ is not integrable on $[0,s]$ with respect to $dx$.
\end{remark}
For every $t\in [0,T]$ the summand $F^{n}_{t}$ can be written as
\begin{eqnarray}
\label{eqdeffn}
F^{n}_{t}&=&\int_{0}^{t}\int_{0}^{t} dB(x_{1}) dB(x_{2}) \int _{(x_{1}\vee x_{2} +\varepsilon _{n})\wedge t } ^{t}ds (s-x_{1} )^{\frac{H}{2}-1}(s-x_{2} )^{\frac{H}{2}-1}\notag\\
&=&\int_{0}^{t-\eps _{n}}\int_{0}^{t-\eps _{n}} dB(x_{1}) dB(x_{2}) \int _{(x_{1}\vee x_{2}) } ^{t-\eps _{n}}ds (s+\varepsilon _{n}-x_{1} )^{\frac{H}{2}-1}(s+\varepsilon _{n}-x_{2} )^{\frac{H}{2}-1}\notag\\
&=& \int_{0}^{t-\eps _{n}} \left( \int_{0} ^{s} (s+\varepsilon _{n}-x) ^{\frac{H}{2}-1} dB(x) \right) ^{2}ds =\int_{0}^{t-\eps _{n}} (Y^{3,H}_{s} )^{2}ds
\end{eqnarray}
where we denoted by
\begin{equation}
\label{y4H-1}Y^{3,H}_{s}= \int_{0} ^{s} (s+\varepsilon _{n}-x) ^{\frac{H}{2}-1} dB(x)\ \ \ \ \text{for} \  0\leq s\leq T.
\end{equation}
 Note that the process $Y^{3,H}$ depends on $n$. But we prefer to use the notation $Y^{3,H}$ without $n$ in order to keep the coherence with the other terms treated before and in the sequel.

Let  $(B_{1}(s)) _{s\in [0,T]}$ the restriction of the Wiener process $B$ to the interval $[0,T]$ and let $Z_{1}^{(n )} $ be the corresponding transport process defined in (\ref{edeftranspro1}) that converges to $B$ in the strong sense (\ref{eq98}). Then,

\begin{equation}
\label{y4H-2}Y^{3,H}_{s}= \int_{0} ^{s} (s+\varepsilon _{n}-x) ^{\frac{H}{2}-1} dB_1(x)
\end{equation}
and  we define
\begin{equation}\label{y4Hn}
Y^{3,H,n}_{s}=  \int_{0} ^{s} (s+\varepsilon _{n}-x) ^{\frac{H}{2}-1} dZ^{(n)}_{1}(x) \ \ \ \ \ \text{for} \ \ 0\leq s\leq t.
\end{equation}
We will show first that $Y^{3,H, n}$ is a strong approximation of $Y^{3,H}$.
\begin{prop}\label{pp3}
Let $Y^{3,H}$, $Y^{3,H,n }$ and $\alpha_n$ be given by (\ref{y4H-2}), (\ref{y4Hn})  and (\ref{eqdefalphan}) respectively,
\begin{equation*}
P \left( \sup_{s\in [0,T] } | Y^{3,H , n}_{s}-Y_{s}^{3,H} | > C \alpha _{n} \right) = o(n^{-q})
\end{equation*}
for each $q>0$ and for $\beta \in (0, \frac{1}{2})$.
\end{prop}
{\bf Proof: } After integrating by parts, we can write, for every $s\in [0,T]$,
\begin{eqnarray*}
| Y^{3,H , n}_{s}-Y^{3,H}_s|&\leq & \eps _{n} ^{\frac{H}{2}-1} \left| B_{1}(s) -Z^{(n)}_{1} (s) \right| \\
&&+ (1-H/2) \int_{0} ^{s} (s+\varepsilon _{n}-x) ^{\frac{H}{2}-2} \left| B_{1}(x)-Z^{(n)}_{1} (x)\right| dx \\
&\leq & 2 \eps _{n}^{\frac{H}{2}-1} \Vert B_{1}-Z^{(n)}_{1} \Vert _{\infty , T}= 2 n ^{\beta} \Vert B_{1}-Z^{(n)}_{1} \Vert _{\infty , T}
\end{eqnarray*}
using the choice of $\eps_{n}$ and hence by (\ref{eq98}),
 \begin{eqnarray*}
P\left( \Vert Y^{3,H , n}-Y^{3,H}\Vert _{\infty, T} \geq C\alpha _{n}\right) \leq P\left(  \Vert B_{1}-Z^{(n)}_{1} \Vert _{\infty , T} \geq C n^{-1/2}(\log n)^{5/2}\right)= o(n^{-q})
\end{eqnarray*}
for every $q>0$ by (\ref{eq98}).
 \qed

\vskip0.2cm

We will introduce now the approximation processes that will converge to $X^{3,H}$. Let us denote, for every $t\in [0,T]$,  by
\begin{equation}
\label{x4n}
X^{3,H,n}_{t}= c(H)\int_{0} ^{t} (Y^{3,H,n} _{s}) ^{2} ds.
\end{equation}

The part $X^{3,H}$ is approximated as follows.
\begin{prop}\label{p4}
For $0<\max \left( \frac{1-\frac{H}{2}}{3-2H}, \frac{2-H}{2+2H}\right) < \beta < \frac{1}{2}$ fixed, let $X^{3,H}$ and $X^{3,H,n}$ defined by (\ref{eqdescx4}) and (\ref{x4n}) respectively. Then for every $\gamma $ such that $0<\gamma<\beta$ and $\gamma+\beta<\frac{1}{2}$,

\begin{equation*}
P \left( \limsup_{n\to \infty} \{ \Vert X^{3,H,n}-X^{3,H} \Vert _{\infty , T} \geq C n^{-(1/2-\beta-\gamma)}(\log n)^{5/2} \} \right) =0.
\end{equation*}
\end{prop}
{\bf Proof: } By (\ref{eqdescx4}) and (\ref{eqdeffn}),
\begin{eqnarray*}
X^{3,H,n}_{t}-X^{3,H}_{t}&=& c(H)\int_{0} ^{t} \left[ (Y^{3,H,n}_{s}) ^{2}-(Y^{3,H}_{s})^{2} \right] ds - c(H)G^{n}_{t} \\
&=& c(H)\int_{0} ^{t} \left[ (Y^{3,H}_{s}-Y^{3,H,n}_{s } ) ^{2} + 2 Y^{3,H}_{s}(Y^{3,H,n}_{s} -Y^{3,H}_s) \right] ds - c(H)G^{n}_{t}
\end{eqnarray*}
we have the bound
\begin{eqnarray*}
\Vert X^{3,H} -X^{3,H,n}\Vert _{\infty ,T} &\leq & c(H)\int_{0}^{T} \vert Y^{3,H}_{s} -Y^{3,H,n}_{s} \vert ^{2}ds \\
&+& 2 c(H)\int_{0}^{T} |Y^{3,H}_{s} | \vert Y^{3,H}_{s} -Y^{3,H,n}_{s} \vert ds + c(H)\sup_{t\in [0,T]} | G^{n}_{t} |.
\end{eqnarray*}

The first two summand in the right hand side above can be treated as in the case of $X^{1,H}$ in Proposition \ref{p2}. We note, in order to apply Lemma  \ref{lims}, we need to notice that $\mathbb{E} \left| \int_{0}^{t} Y^{3,H}_{s} ds \right| < C$ with $C$ not depending on $n$ (Lemma \ref{lims} can still be used although the process $Y^{3,H}$ depends on $n$). Let us handle the term $G^{n}$.  We will actually show that
\begin{equation}
\label{bn}
\sum_{n} P( \sup_{t\in [0,T]} G^{n}_{t}>C n^{-(1/2-\beta-\gamma)}(\log n)^{5/2})<\infty.
\end{equation}
which will imply that
\begin{equation*}
P \left( \limsup_{n}\{ \Vert G^{n}\Vert _{\infty, T} > C  n^{-(1/2-\beta-\gamma)}(\log n)^{5/2}\} \right) =0.
\end{equation*}
For every $t\in [0,T]$ we have, using the fact that the integrand is symmetric with respect to the variables $x_{1}$ and $x_{2}$,
\begin{eqnarray*}
G^{n}_{t}&=& 2\int_{0}^{t} dB(x_{1}) \int_{0}^{x_{1} } dB(x_{2}) \int_{x_{1}} ^{(x_{1}+ \varepsilon _{n}) \wedge t} ds (s-x_{1}) ^{\frac{H}{2}-1} (s-x_{2}) ^{\frac{H}{2}-1}\\
&=& 2 \int_{0} ^{(t-\varepsilon _{n})\vee 0} dB(x_{1}) \int_{0} ^{x_{1}} dB(x_{2}) \int_{x_{1}} ^{ x_{1}+\eps _{n} } ds (s-x_{1}) ^{\frac{H}{2}-1} (s-x_{2}) ^{\frac{H}{2}-1}\\
&&+ 2 \int_{(t-\eps _{n})\vee0} ^{t} dB(x_{1}) \int _{0}^{x_{1}} dB(x_{2})  \int_{x_{1}} ^{t} ds  (s-x_{1}) ^{\frac{H}{2}-1} (s-x_{2}) ^{\frac{H}{2}-1}\\
&:=& G^{1,n}_{t} + G^{2,n}_{t} .
\end{eqnarray*}
Note that the mapping
$$x_{1} \to \int _{0}^{x_{1}} dB(x_{2})  \int_{x_{1}} ^{ x_{1}+\eps _{n} } ds  (s-x_{1}) ^{\frac{H}{2}-1} (s-x_{2}) ^{\frac{H}{2}-1}$$
is adapted with respect to ${\cal{F}} _{x_{1}} $ (the filtration generated by the Wiener process $B$). Then the process $(G^{1,n} _{t}) _{t}$ is a martingale  for every $n$. Then we have, taking $\hat{\alpha}_n=n^{-(1/2-\beta-\gamma)}(\log n)^{5/2}$ and using Doob's inequality,

\begin{eqnarray*}
&&P \left( \sup_{t\in [0,T] } G^{1,n}_{t} \geq C n^{-(1/2-\beta-\gamma)}(\log n)^{5/2}\right) \leq \hat{\alpha}_n ^{-p} \mathbb{E} \sup _{t\in [0,T]}\left| G^{1,n }_{t} \right| ^{p} \\
&\leq & C \hat{\alpha}_n^{-p} \mathbb{E} \left|\int_{0} ^{(T-\varepsilon)\vee 0} dB(x_{1}) \int_{0} ^{x_{1}}dB(x_{2}) \int_{x_{1}}^{x_{1}+\eps _{n}}  ds(s-x_{1}) ^{\frac{H}{2}-1} (s-x_{2}) ^{\frac{H}{2}-1}\right| ^{p}
\end{eqnarray*}
with $C$ allowed to depend also on $p$ in this proof. Note that the random variable \\
 $ \int_{0} ^{(T -\eps_{n})\vee 0} dB(x_{1}) \left( \int_{0} ^{x_{1}}dB(x_{2}) \int_{x_{1}}^{x_{1}+\eps _{n}}  ds(s-x_{1}) ^{\frac{H}{2}-1} (s-x_{2}) ^{\frac{H}{2}-1}\right) ^{2}$ is a multiple integral of order two. Therefore, by the hypercontractivity property (\ref{hyper}), it is not difficult to see that

\begin{eqnarray*}
&&P \left( \sup_{t\in [0,T] } G^{1,n}_{t} \geq C n^{-(1/2-\beta-\gamma)}(\log n)^{5/2}\right) \\
 &&\leq C \hat{\alpha}_n ^{-p}
 \left[ \int_{0}^{(T -\eps_{n})\vee 0} dx_{1}   \int_{0}^{x_{1}} dx_{2} \left( \int_{x_{1}}^{x_{1}+\eps _{n}}  ds(s-x_{1}) ^{\frac{H}{2}-1} (s-x_{2}) ^{\frac{H}{2}-1}\right) ^{2}\right] ^{\frac{p}{2}}
\end{eqnarray*}
Let us first compute the integral with respect to $ds$. By making the change of variables $z=\frac{s-x_{1}}{s-x_{2}}$ with $ds = \frac{x_{1}-x_{2}}{(1-z)^{2}}dz$ we get

\begin{eqnarray*}
&&\int_{0}^{T-\eps _{n}
} dx_{1}   \int_{0}^{x_{1}} dx_{2} \left( \int_{x_{1}}^{x_{1}+\eps _{n}}  ds(s-x_{1}) ^{\frac{H}{2}-1} (s-x_{2}) ^{\frac{H}{2}-1}\right) ^{2}\\
&=& \int_{0}^{T-\eps _{n}
} dx_{1}   \int_{0}^{x_{1}} dx_{2}  (x_{1}-x_{2}) ^{2H-2} \left(  \int_{0} ^{\frac{\eps _{n}}{\eps _{n} + x_{1}-x_{2}}}z^{\frac{H}{2}-1} (1-z) ^{-H}dz\right) ^{2}.
\end{eqnarray*}
We separate the integral $dx_{1}dx_{2}$ into two regions: when $x_{1}-x_{2} \leq\eps _{n}$ and when $x_{1}-x_{2} > \eps _{n}$. The above term will be bounded by

\begin{eqnarray*}
&& \int_{0}^{T} dx_{1}   \int_{(x_{1}-\eps _{n})\vee 0}^{x_{1}} dx_{2}  (x_{1}-x_{2}) ^{2H-2}
 \left(  \int_{0} ^{\frac{\eps _{n}}{\eps _{n} + x_{1}-x_{2}}}z^{\frac{H}{2}-1} (1-z) ^{-H}dz\right) ^{2}\\
 &&+ \int_{0}^{T} dx_{1}   \int_{0}^{(x_{1}-\eps _{n})\vee 0} dx_{2}  (x_{1}-x_{2}) ^{2H-2}
 \left(  \int_{0} ^{\frac{\eps _{n}}{\eps _{n} + x_{1}-x_{2}}}z^{\frac{H}{2}-1} (1-z) ^{-H}dz\right) ^{2}\\
 &\leq &  \int_{0}^{T} dx_{1}   \int_{(x_{1}-\eps _{n})\vee 0}^{x_{1}} dx_{2}  (x_{1}-x_{2}) ^{2H-2}
 \left(  \int_{0} ^{1}z^{\frac{H}{2}-1} (1-z) ^{-H}dz\right) ^{2}\\
 &&+  \int_{0}^{T} dx_{1}   \int^{(x_{1}-\eps _{n})\vee 0}_{0} dx_{2}  (x_{1}-x_{2}) ^{2H-2} \left(  \int_{0} ^{\frac{\eps _{n}}{\eps _{n} + x_{1}-x_{2}}}z^{\frac{H}{2}-1} (1-z) ^{-H}dz\right) ^{2}\\
 &\leq&c(H)  \int_{0}^{T} dx_{1}   \int_{x_{1}-\eps _{n}}^{x_{1}} dx_{2}  (x_{1}-x_{2}) ^{2H-2}\\
 &&+ c(H) \int_{0}^{T} dx_{1}   \int^{(x_{1}-\eps _{n})\vee 0}_0 dx_{2}  (x_{1}-x_{2}) ^{2H-2} F_1^2(H/2, H,  H/2+1, 1)\left(\frac{\eps _{n}}{\eps _{n} + x_{1}-x_{2}}\right) ^{H} \\
 &\leq & C \eps _{n} ^{2H-1},
\end{eqnarray*}
where $F_1^2(H/2, H,  H/2+1, 1)$ is the incomplete beta function and hence
\begin{eqnarray*}
&&P \left( \sup_{t\in [0,T] } G^{1,n}_{t} \geq c\hat{\alpha} _{n}  \right) \leq C \hat{\alpha} _{n} ^{-p} \eps _{n} ^{(2H-1)\frac{p}{2}}
\end{eqnarray*}
and the series $\sum_{n} \eps _{n} ^{\frac{p}{2} (2H-1) } \hat{\alpha} _{n} ^{-p}$ is finite if
$$p \left( \beta \frac{2H-1}{2-H} -(\frac{1}{2}-\beta -\gamma ) \right) >1.$$
Note that $\beta >\frac{2-H}{2+2H}$ implies that $ \beta \frac{2H-1}{2-H} -(\frac{1}{2}-\beta -\gamma >0$ for small $\gamma >0$. By choosing $p$ large enough we obtain that
\begin{equation*}
\sum_{n}P \left( \sup_{t\in [0,T] } G^{1,n}_{t} \geq C\hat{\alpha}_n \right) <\infty
\end{equation*} for every $\beta \in (0, \frac{1}{2})$.

Let us handle now the term denoted by $G^{2,n}$. We have
\begin{eqnarray*}
P \left( \sup_{t\in [0,T] } G^{2,n}_{t} \geq C\hat{\alpha} _{n}  \right)&\leq & \hat{\alpha} _{n} ^{-p} \mathbb{E} \sup _{t\in [0,T]}\left| G^{2,n }_{t} \right| ^{p} \\
\end{eqnarray*}
In order to control $\mathbb{E} \sup _{t\in [0,T]}\left| G^{2,n }_{t} \right| ^{2}$ we will use Garsia's lemma. To this end we need to estimate the $L^{p}$ norm of the increment $G^{2,n}_{t}-G^{2,n}_{s}$ when $t$ is close to $s$. Note first that, by the change of variables $z=\frac{s-x_{1}}{s-x_{2}}$ we have
\begin{equation*}
G^{2,n}_{t}= \int_{t-\eps _{n}} ^{t} dB(x_{1})\int_{0}^{x_{1}}dB(x_{2}) \vert x_{1}-x_{2} \vert ^{H-1} \int_{0} ^{\frac{t-x_{1}}{t-x_{2}}} z^{\frac{H}{2}-1}(1-z) ^{ -H} dz
\end{equation*}
and for $t,s \in [0,T]$ such that $s>t-\eps _{n}$, by the isometry of multiple stochastic integrals (\ref{eqisometry})
\begin{eqnarray*}
&&\mathbb{E}\left|G^{2,n}_{t}-G^{2,n}_{s}\right| ^{2}\\
&=&c(H) \int_{t-\eps_{n}}^{t}dx_{1} \int_{0} ^{x_{1}} dx_{2} \vert x_{1}-x_{2}\vert ^{2H-2} \left(\int_{0} ^{\frac{t-x_{1}}{t-x_{2}}} z^{\frac{H}{2}-1}(1-z) ^{ -H} dz \right)^{2} \\
&&+ c(H) \int_{s-\eps_{n}}^{s}dx_{1} \int_{0} ^{x_{1}} dx_{2} \vert x_{1}-x_{2}\vert ^{2H-2} \left(\int_{0} ^{\frac{s-x_{1}}{s-x_{2}}} z^{\frac{H}{2}-1}(1-z) ^{ -H} dz \right)^{2}\\
&&-2c(H) \int_{t-\eps}^{s} dx_{1} \int_{0} ^{x_{1}}dx_{2}  \vert x_{1}-x_{2}\vert ^{2H-2}
\left(\int_{0} ^{\frac{s-x_{1}}{s-x_{2}}} z^{\frac{H}{2}-1}(1-z) ^{ -H} dz \right)\left(\int_{0} ^{\frac{t-x_{1}}{t-x_{2}}} z^{\frac{H}{2}-1}(1-z) ^{ -H} dz \right)
\end{eqnarray*}
and by majorizing $\int_{0} ^{\frac{t-x_{1}}{t-x_{2}}} z^{\frac{H}{2}-1}(1-z) ^{ -H} dz$ by $\int_{0} ^{1} z^{\frac{H}{2}-1}(1-z) ^{ -H} dz= \beta (\frac{H}{2}, 1-H)$ we obtain the bound
\begin{eqnarray*}
\mathbb{E}\left|G^{2,n}_{t}-G^{2,n}_{s}\right| ^{2}
&\leq&c(H) \int_{t-\eps_{n}}^{t}dx_{1} \int_{0} ^{x_{1}} dx_{2} \vert x_{1}-x_{2}\vert ^{2H-2}\\
&&+ c(H) \int_{s-\eps_{n}}^{t}dx_{1} \int_{0} ^{x_{1}} dx_{2} \vert x_{1}-x_{2}\vert ^{2H-2}+ c(H) \int_{s-\eps_{n}}^{t}dx_{1} \int_{0} ^{x_{1}} dx_{2} \vert x_{1}-x_{2}\vert ^{2H-2}\eps _{n}.
\end{eqnarray*}
The hypercontractivity property of multiple integrals (\ref{hyper}) implies that for all $t,s$
\begin{eqnarray*}
\mathbb{E}\left|G^{2,n}_{t}-G^{2,n}_{s}\right| ^{p}\leq c(p,H) \left( \mathbb{E}\left|G^{2,n}_{t}-G^{2,n}_{s}\right| ^{2}\right) ^{\frac{p}{2} } \leq C \eps _{n} ^{\frac{p}{2}}.
\end{eqnarray*}
where $C$ is a constant that depend on $
H$ and $p$. Finally, by Garsia lemma (see e.g. \cite{N}, Appendix A.3) for $p>2$
\begin{equation}\label{fin1}
\mathbb{E} \sup_{0\leq t\leq T} \vert G^{2,n}\vert ^{p} \leq C\eps _{n} ^{\gamma}
\end{equation}
for every $\gamma $ such that $0<\gamma < \frac{p}{2}-1$.   The bound (\ref{fin1}) implies, using Markov's inequality and taking suitable $p$ large enough,  that
\begin{equation}\label{nnn}
\sum_{n}P \left( \sup_{t\in [0,T] } G^{2,n}_{t} \geq c\hat{\alpha _{n} } \right) <\infty
\end{equation}
due to the fact that $\beta > \frac{1-\frac{H}{2}}{3-2H}$ and this finishes the proof. \qed

\vskip0.3cm

Let us finally treat the summand $X^{2,H}$ in (\ref{defx1x4}).  Its approximation will be a mixture of the approximations of  $X^{1,H}$ and $ X^{3,H}$. We have
\begin{eqnarray*}
X^{2,H}_{t}&=&\int_{-\infty}^{0} dB(x_{1}) \int_{0} ^{t} dB(x_{2}) \int_{x_{2} }^{t} (s-x_{1}) ^{\frac{H}{2}-1} (s-x_{2}) ^{\frac{H}{2}-1} ds \\
&=& \int_{0} ^{t} dB(x_{2}) \int_{x_{2}} ^{t} Y^{1,H}_{s} (s-x_{2} ) ^{\frac{H}{2}-1} ds \\
\end{eqnarray*}
with $Y^{1,H}$ given by (\ref{y1H}). To avoid the singularity of the integral with respect to $ds$ at $s=x_{2}$ we will decompose this integral into two parts. In this way we can write, for $\eps _{n}$ the sequence converging to 0 as $n\to \infty$ chosen  before
 \begin{eqnarray}
 \label{eqdefx2}
&&X^{2,H}_{t}= \int_{0} ^{t} dB(x_{2}) \int_{(x_{2}+\eps _{n})\wedge t} ^{t} Y^{1,H}_{s} (s-x_{2} ) ^{\frac{H}{2}-1} ds \notag \\
&&+ \int_{0} ^{t} dB(x_{2}) \int_{x_{2}} ^{(x_{2}+ \eps _{n} )\wedge t} Y^{1,H}_{s} (s-x_{2} ) ^{\frac{H}{2}-1} ds \notag\\
&=& \int_{\eps _{n}}^{t} ds Y^{1,H}_{s}  \int_{0} ^{s-\eps _{n}} (s-x_{2}) ^{\frac{H}{2}-1} dB(x_{2}) +\int_{(t-\eps _{n})\vee 0 } ^{t} dB(x_{2}) \int_{x_{2}}^{(x_2 + \eps _{n})\wedge t}  Y^{1,H}_{s} (s-x_{2} ) ^{\frac{H}{2}-1} ds\notag\\
&&+\int^{(t-\eps _{n})\vee 0 } _{0} dB(x_{2}) \int_{x_{2}}^{(x_2 + \eps _{n})\wedge t}  Y^{1,H}_{s} (s-x_{2} ) ^{\frac{H}{2}-1} ds\notag\\
&=& \int_{\eps _{n}} ^{t} Y^{1,H}_{s} Y^{', 3,H}_{s} ds + \mathbf{F}^{n}+\mathbf{G}^{n}
\end{eqnarray}
with
\begin{equation}\label{yprime}
Y^{', 3,H}_{s} = \int_{0} ^{s-\eps _{n}} (s-x_{2}) ^{\frac{H}{2}-1} dB(x_{2}),  Y^{', 3,H,n }_{s} = \int_{0} ^{s-\eps _{n}} (s-x_{2}) ^{\frac{H}{2}-1}dZ^{(1)}(s), \hskip0.5cm s\in [0,T]
 \end{equation}
As in the proof of (\ref{bn}) and (\ref{nnn}] we can show that
\begin{equation*}
\sum_{n}P \left( \sup_{t\in [0,T] } \mathbf{F}^{n}_{t} \geq c\hat{\alpha} _{n} ^{2} \right) <\infty \mbox{ and } \sum_{n}P \left( \sup_{t\in [0,T] } \mathbf{G}^{n}_{t} \geq c\hat{\alpha} _{n} ^{2} \right) <\infty
\end{equation*}
The approximation result  to $X^{2,H}$ is stated in the next proposition. We will use the process $Y^{3,H, n}$ instead of $Y^{',3,H,n}$ because clearly they are very close and one can replace the other.
\begin{prop}\label{pp5}
For $0< \max \left( \frac{1-H/2}{3-2H}, \frac{2-H}{2+2H}\right)< \beta < \frac{1}{2}$ fixed, let  $Y^{1,H,n}$, $Y^{3,H,n}$ and $X^{2,H}$  be given by (\ref{eqdefYn}), (\ref{y4Hn}) and (\ref{eqdefx2}) respectively. Define
\begin{equation}
\label{x2Hn}
X^{2,H,n}_{t}= \int_{0}^{t} ds Y^{1,H,n}_{s} Y^{3,H,n}_{s} , \hskip0.5cm t\in [0,T].
\end{equation}
Then for every $\gamma $ such that $0<\gamma<\beta$ and $\gamma+\beta<\frac{1}{2}$,

\begin{equation*}
P\left( \limsup_{n} \{ \Vert X^{2,H,n}-X^{2,H} \Vert _{\infty ,T}>Cn^{-(1/2-\beta-\gamma)}(\log n)^{5/2} \} \right) =0.
\end{equation*}
\end{prop}
{\bf Proof: }The proof follows from the proofs of Proposition \ref{p2} and Proposition \ref{p4} since for every $s$ we have $$2 Y^{1,H,n}_{s} Y^{', 3,H,n}_{s}\leq (Y_{s}^{1,H,n})^{2}+ (Y_{s}^{', 3,H,n}) ^{2}.$$\qed

\vskip0.2cm

Let us summarize the conclusions of Proposition \ref{p2}, \ref{p4} and \ref{pp5} in the main result of our paper.

\begin{theorem}
Let $X^{H} $ be the Rosenblatt process (\ref{rose}) and $0<\max \left( \frac{1-H/2}{3-2H}, \frac{2-H}{2+2H}\right) < \beta < \frac{1}{2}$ fixed. Define
\begin{equation*}
X^{H,n}_{t}= X^{1,H,n}_{t} + 2X^{2,H,n}_{t} + X^{3,H,n}_{t}, \hskip0.5cm t\in [0,T]
\end{equation*}
with $X^{1,H,n}, X^{2,H,n}, X^{3,H,n}$ given by (\ref{x1n}), (\ref{x2Hn}), (\ref{x4n}) respectively. Then for every $\gamma $ such that $0<\gamma<\beta$ and $\gamma+\beta<\frac{1}{2}$,
\begin{equation*}
P\left( \limsup_{n} \{ \Vert X^{H,n}-X^{H} \Vert _{\infty ,T}>Cn^{-(1/2-\beta-\gamma)}(\log n)^{5/2}  \} \right) =0.
\end{equation*}
\end{theorem}

\begin{remark}
The slowest rate of convergence is obtained for $H$ close to one because in this case $\beta $ is close to $\frac{1}{2}$. When $H$ is close to $\frac{1}{2}$ then $\beta $ is close to $\frac{3}{8}$. But this situation cannot be compared with previous results in the literature because the Rosenblatt process is not defined for $H=\frac{1}{2}$.
\end{remark}


\begin{thebibliography}{99}
\bibitem {BN}{J.-C. Breton and I. Nourdin (2008): }Error bounds on the
non-normal approximation of Hermite power variations of fractional Brownian
motion. {\emph{Electronic Communications in Probability}}, \textbf{13}, 482-493.

\bibitem {BrMa}{P. Breuer and P. Major (1983): }Central limit theorems for
nonlinear functionals of Gaussian fields. \emph{J. Multivariate Analysis}{,
\textbf{13 }(3), 425-441.}

\bibitem {CNT}{A. Chronopoulou, C.A. Tudor and F. Viens (2009): }Application
of Malliavin calculus to long-memory parameter estimation for non-Gaussian
processes. \emph{Comptes rendus - Mathematique} \textbf{347}, 663-666.

\bibitem{CH}
{M. Cs\"orgo and L. Horvath (1988): }{\em Rate of convergence of transport processes with an application to stochastic differential equations. } Probability Theory and Related Fields, {\bf 78}, 379-387.

\bibitem {DM}{R.L. Dobrushin and P. Major (1979): }\emph{Non-central limit
theorems for non-linear functionals of Gaussian fields. }{Z.
Wahrscheinlichkeitstheorie verw. Gebiete, \textbf{50}, 27-52. }

\bibitem{E}
{N. Enriquez (2004): }{\em A simple construction of the fractional Brownian  motion. } Stochastic Processes and their Applications, {\bf 109}, 203-223.

\bibitem{GaGoLe}
{J. Garzon, L.G. Gorostiza and J.A. Leon (2009): }{\em A strong uniform approximation of fractional Brownian motion by means of  transport processes. } Stochastic Processes and their Applications, {\bf 119}, 3435-3452.

\bibitem{GoGr}
{L.G. Gorostiza, R.J. Griego (1980): }{\em Rate of convergence of uniform transport processes to Brownian motion and applications to stochastic integrals. }Stochastics, {\bf 3}, 291-303.

\bibitem{GHR}
{R.J. Griego, D. Heath, A. Ruiz-Moncayo (1971): }{\em Almost sure convergence of uniform transport processes  to Brownian motion. } Ann. Math. Stat. {\bf 42}, 1129-1131.

\bibitem{Hall}
{P. Hall, W. Hardle, T. Kleinow and P. Schmidt (2000): } Semiparametric Bootstrap Approach to Hypothesis tests and Confidence
intervals for the Hurst coefficient. \emph{Stat. Infer. Stoch. Process.} \textbf{3},   263–276.

\bibitem{LK} {A.J. Lawrance and N.T. Kottegoda (1977): }Stochastic modelling of riverflow time series. \emph{J. Roy. Statist. Soc. Ser. A }, \textbf{140}(1), 1-47.


      \bibitem{Major} {\sc Major, P.} (2005). {\em Tail behavior of multiple
    random integrals and $U$-statistics. }{Probability Surveys.}

\bibitem{N} {D. Nualart (2006): }\emph{Malliavin Calculus and Related
Topics. Second Edition. }{Springer. }

\bibitem {NNT}{I. Nourdin, D. Nualart and C.A Tudor (2007): }\emph{Central and
Non-Central Limit Theorems for weighted power variations of the fractional
Brownian motion. }  Annales I.H.P. -Probabilit\'es et Statistiques, {\bf 46} (4), 1055-1079.

\bibitem{Sz}
{T. Szabados (1996): }{\em Strong approaximation of fractional Brownian motion by movin averages of simple random walks. } Stochastic processes and their applications, {\bf 31}, 243-255.

\bibitem {Ta1}{M. Taqqu (1975): }Weak convergence to the fractional
Brownian motion and to the Rosenblatt process. \emph{Z. Wahrscheinlichkeitstheorie
verw. Gebiete,}{ \textbf{31}, 287-302.}

\bibitem{Taqqu3} {M. Taqqu (1978): } A representation for self-similar processes. \emph{ Stochastic Processes and their Applications, }\textbf{7}, 55-64.

\bibitem {T}{C.A. Tudor (2008): }Analysis of the Rosenblatt process.
\emph{ESAIM Probability and Statistics}{, \textbf{12}, 230-257.}

\bibitem {TV}{C.A. Tudor and F. Viens (2008}): Variations and estimators
through Malliavin calculus. Annals of Probability, {\bf 37} (6), 2093-2134.

\bibitem{Wu}{W.B. Wu (2005): } Unit root testing for functionals of linear processes. \emph{Econ. Theory, } \textbf{22},  1–14.
\end{thebibliography}
\end{document}